\documentstyle[12pt,amssymb]{amsart}
\input{diagrams.tex}
\diagramstyle[scriptlabels,height=8mm,width=8mm]
\theoremstyle{plain}

\newcommand{\query}[1]%
{\mbox{}\marginpar{\raggedright\hspace{0pt}{\small\em #1}}}%

\newtheorem{thm}{Theorem}[section]
\newtheorem{cor}[thm]{Corollary}
\newtheorem{lem}[thm]{Lemma}
\newtheorem{prop}[thm]{Proposition}

\theoremstyle{definition}
\newtheorem{defi}[thm]{Definition}
\newtheorem{conj}[thm]{Conjecture}
\newtheorem{conv}[thm]{Convention}
\newtheorem{nota}[thm]{Notation}
\newtheorem{rem}[thm]{Remark}
\newtheorem{exa}[thm]{Example}
\newtheorem{sit}[thm]{}
\newcommand{\new}{\newcommand}

\new{\f}[1]{\frak {#1}}
\def\cL{{\cal L}}
\def\cO{{\cal O}}

\new{\dt}{,\dots,}
\new{\si}{\sigma}
\new{\Si}{\Sigma}
\new{\h}[1]{\widehat {#1}}
\new{\bari}[1]{{#1}^{(i)}}
\new{\barun}[1]{{#1}^{(1)}}
\new{\barn}[1]{{#1}^{(n)}}
\new{\quot}[2]{\raise 2pt\hbox{#1}/\raise-4pt\hbox{#2}}
\new{\xb}{\bar x}
\new{\yb}{\bar y}
\new{\zb}{\bar z}

\new{\ga}{\gamma}

\newcommand{\pP}{{\mathbb P}}
\newcommand{\A}{{\mathbb A}}

\newcommand{\C}{{\mathbb C}}
\newcommand{\Q}{{\mathbb Q}}
\newcommand{\Z}{{\mathbb Z}}

\newcommand{\G}{{\Gamma}}

\newcommand{\caO}{{\cal {O}}}

\newcommand{\pr}{{\bf P}}
\newcommand{\Aut}{{\operatorname{Aut}\,}}

\newcommand{\baf}{{\bar {f}}}

\newcommand{\proof} {\no{\it Proof. ~} }
\newcommand{\brem}{\begin{rem}}
\newcommand{\erem}{\end{rem}}
\newcommand{\bexa}{\begin{exa}}
\newcommand{\eexa}{\end{exa}}
\newcommand{\bdefi}{\begin{defi}}
\newcommand{\edefi}{\end{defi}}
\newcommand{\bcor}{\begin{cor}}
\newcommand{\ecor}{\end{cor}}
\newcommand{\blem}{\begin{lem}}
\newcommand{\elem}{\end{lem}}
\newcommand{\bconv}{\begin{conv}}
\newcommand{\econv}{\end{conv}}
\newcommand{\bconj}{\begin{conj}}
\newcommand{\econj}{\end{conj}}
\newcommand{\bprop}{\begin{prop}}
\newcommand{\eprop}{\end{prop}}
\newcommand{\bthm}{\begin{thm}}
\newcommand{\ethm}{\end{thm}}
\newcommand{\bnota}{\begin{nota}}
\newcommand{\enota}{\end{nota}}
\newcommand{\bsit}{\begin{sit}}
\newcommand{\esit}{\end{sit}}
\newcommand{\be}{\begin{eqnarray}}
\newcommand{\ee}{\end{eqnarray}}
\newcommand{\no}{\noindent}

\edef\qqed{\qed}
\def\qed{\qqed \medskip}

\begin{document}

\title{Families of affine planes:\\
the existence of a cylinder}

\vspace{.5cm}

\author{Shulim Kaliman}
\address{Department of Mathematics,
University of Miami,
Coral Gables, FL  33124, U.S.A.} 
\email{kaliman@@math.miami.edu}

\author{Mikhail Zaidenberg}
\address{Universit\'e
Grenoble I, Institut Fourier, UMR 5582 CNRS-UJF, BP 74,
38402 St.\ Martin
d'H\`eres c\'edex, France}
\email{zaidenbe@@ujf-grenoble.fr}

\thanks{Research of the  first author are
partially supported by the NSA grant  MDA904-00-1-0016.\\
This work was finished during the second author's stay at the Max Planck
Institute of Mathematics at Bonn; 
he would like to thank MPI for the hospitality. \\
\mbox{\hspace{11pt}}{\it 2000 Mathematics Subject Classification}: 
14R25, 14R05, 14D06, 14E99.\\
\mbox{\hspace{11pt}}{\it Key words}: affine plane, Ramanujam-Morrow graph,
families of quasi-projective varieties, contraction theorem} 

\date{}

\begin{abstract} Given a family of complex affine planes, we 
show that it is trivial over 
a Zariski open subset of the base. 
The proof relies upon a relative version
of the contraction theorem. 
\end{abstract}

\maketitle

\medskip

\begin{center}

{\large Introduction}

\end{center}
\medskip

The varieties $X,S,U,V$ considered below are 
assumed being smooth, quasi-projective and defined over $\C$.
By a {\em cylinder over} $U$ we mean a Cartesian 
product $U \times \C^k$
where $k >0$. We say that a family $f : X \to S$ of
quasi-projective varieties {\em contains
a cylinder} if, for some Zariski open
subset $S_0$ of $S$, there is a commutative diagram  

\begin{diagram}[w=7mm,h=8mm,midshaft]
f^{-1}(S_0)
&& \rTo^{\varphi} &&  S_0 \times \C^k\\
&\rdTo<f &&
\ldTo>{pr_1}\\
&&  S_0 \,
\end{diagram}

\noindent where $\varphi$ is an isomorphism. 

The main result of this paper is the following 

\no\bthm\label{mt} A smooth family $f : X \to S$ with
general fibers isomorphic to $\C^2$ contains a cylinder 
$S_0 \times \C^2$.\ethm

A similar result is well-known
for families of relative dimension 1 
\cite{BhaDut, EakHei, KamMiy, 
KamWri, Miy1, Miy2, MiySug, Sug}.
In relative dimension 2 the existence 
of a cylinder was   
an assumptions in the previous work of 
\cite{AsaBha, Kam, Miy2, Sat}. 
For higher relative dimensions, 
it makes a part of Conjecture 3.8.3 in 
\cite{DolVei}.

Theorem \ref{mt} proves Lemma III of \cite{Kal} 
which is one of the principal
ingredients in the proof of the following result.
\medskip

\no {\bf Theorem} [Kal]. 
{\em A polynomial $p$ on $\C^3$ with general fibers
isomorphic to $\C^2$
is a variable of the polynomial algebra $\C^{[3]}$.
In particular, all its fibers are isomorphic to $\C^2$.}
\medskip

Section 1 is devoted to a relative version 
of the contraction theorem
over an affine base, which we need in sections 2 and 3
where we prove Theorem \ref{mt}.

We are grateful to 
L. Bonavero, I. Dolgachev, H. Flenner, Sh. Ishii, V. Lin, D. Markushevich  
and P. Russell 
for
their advice and references, 
especially helpful in section 1 below.

\section{A contraction theorem}

The main result of this section (Theorem \ref{pr1}) is a
relative version of the classical 
Castelnuovo-Enriques-Kodaira contraction theorem.
In the analytic setting it follows from
the Moishezon-Nakano-Fujiki theorem 
(see \cite{AncTom, Kod, Moi, FujNak, Nak}
and especially \cite[Rem. 3]{Fuj}), whereas in the projective
setting it follows from the theorem on contraction of extremal 
rays as given in
\cite[Thm. 3-2-1]{KMM} (cf. also \cite{Art, Maz, Cor, Ish}). 
Actually, the particular version 
that we need is much simpler, so we provide a proof along  
the lines of the Castelnuovo-Enriques-Kodaira contraction
theorem \cite[Appendix]{Kod}, \cite[Sect. 4.1, p. 154]{GriHar}.

As usual, for an algebraic variety $X$ its structure sheaf
is denoted by $\cO_X$. If $L$ is a line bundle
on $X$ and $Y$ is a subvariety of $X$ then $\cO_Y(L)$ 
denotes the sheaf of germs of section of $L$ over $Y$.
We begin with the following lemma 
(cf. \cite[\S7.6]{Gra}, \cite[Thm. 4.7]{Pet}).

\medskip\no\blem\label{l1} Let $\rho : E \to S$ be a smooth proper
morphism of smooth quasi-projective varieties with fibers $E_s : = 
\rho^{-1} (s) \, (s \in S)$.
For a line bundle $L$ on $E$ we denote $\cL : =\cO_E (L)$
and $\cL_s :=\cO_{E_s}(L)$. Suppose that
$$H^q(E_s,\cL_s)=0 \qquad \forall s \in S, \quad\forall q \geq 1. 
\leqno {\rm (o)}$$
Then for any Zariski open affine subset $S_0 \subset S$ 
and $E_{S_0} : = \rho^{-1}(S_0)$ we have
$$H^q(E_{S_0},\cL )=0 \qquad \forall q \geq 1\,,
\leqno {\rm (a)}$$         
and

\noindent {\rm (b)} for every point $s \in S_0$ 
the restriction homomorphism
$$H^0(E_{S_0}, \cL ) \to H^0(E_s, \cL_s )$$ 
is surjective. Furthermore,

\noindent {\rm (c)} the sheaf $\rho_* \cL$ 
is locally free and
generated by a vector bundle, say, $\xi(L)$ over $S$
with fibers $\xi(L)_s=H^0(E_s, \cL_s ),\,\,\,s\in S$.
\elem

\no\proof (a) Note that by \cite[Prop. III.9.2.c]{Ha} $\cL$ is a flat
$\cO_S$-module. In virtue of the
assumption (o), for every $s \in S$ and every $q \geq 1$ the
natural homomorphism 
$$R^q\rho_*\cL \otimes_{\cO_S} k(s) \to H^q (E_s, \cL_s)=0$$
is an isomorphism, and the coherent sheaf $R^q\rho_*\cL$ is locally free  
\cite[Thm. 4]{Gra}, \cite[Thm. III.12.11.a, Ex. II.5.8.c]{Ha}, 
\cite[Prop. II.3.7]{Dan} (here $k(s)\simeq \C$ denotes the residue
field of a closed point $s \in S$). 
Thus we have 
$$ R^q\rho_*\cL=0\qquad\forall q\ge 1\,.$$
The Leray spectral sequence gives now isomorphisms
\be\label{1} H^q(E, \cL ) \simeq H^q (S, \rho_* \cL ) \qquad 
\forall q \geq 1 \ee
\cite[(5.16)]{Pet}.
For a Zariski open affine subset $S_0 \subset S$,
by Serre's vanishing theorem 
\cite[Thm. III.3.7]{Ha} we have
$$H^q (S_0, \rho_* \cL )=0 \qquad \forall q \geq 1,$$
which together with (\ref{1})
implies (a).

\no (b) Since $R^1\rho_* \cL =0$, 
for every point $s \in S_0$ 
the homomorphism
\be\label{2} \rho_* \cL \otimes_{\cO_{S_0}}k(s) \to H^0 (E_s, \cL_s )
\ee
is surjective, whence it is an isomorphism \cite[Thm. III.12.11]{Ha}.
On the other hand, since $S_0$
is affine we have an isomorphism 
\be\label{3} 
(\rho_* \cL )|_{S_0} \simeq H^0(E_{S_0}, \cL )^{\sim } \ee
where $M^{\sim }$
denotes the $\cO_{S_0}$-module generated by an $H^0(S_0,
\cO_{S_0})$-module $M$ \cite[Prop. III.8.5]{Ha}. 
Now (\ref{2}) and (\ref{3}) yield 
(b). 

\no (c) By (o) we have $h(s):=\dim H^0(E_s, \cL_s )=\chi(E_{s}, \cL_s )$
where the Euler characteristic
is locally constant on $S$ \cite[Thm. 5]{Gra},
\cite[Prop. II.3.8]{Dan}. Now the isomorphism in (2) and
\cite[Ex. II.5.8.c]{Ha} imply
that $\rho_* \cL$ is a locally free sheaf 
with $\rho_* \cL\simeq \cO_S(\xi(L))$
where $\xi(L)_s=H^0(E_s, \cL_s ),\,\,\,s\in S$. 
The proof is completed. \qed

\no\bcor\label{cor1} 
The statements of Lemma 1.1 remain true if one replaces
the assumption ${\rm (o)}$ by any one of the following two:

\noindent $\rm (o')$ For every $s \in S$ the line bundle 
$L|_{E_s} -K_{E_s}$ on the fiber $E_s$ is ample,
where $K_{E_s}$ is the canonical bundle on $E_s$.

\noindent $\rm (o'')$ For every $s \in S$ and for
each irreducible component $C$ of
the fiber $E_s$, $C\simeq \pP^n \, (n\geq 1)$ and
$L|_C \simeq \cO (l)$ with $l \geq 0$.
\ecor

Indeed, by the Kodaira-Nakano vanishing theorem any one
of the conditions $\rm (o')$ and $\rm (o'')$ 
implies (o) \cite[Sect.1.2, p. 154]{GriHar}.

\medskip \no\bthm\label{pr1} 
Let $\pi : V \to S$ be a smooth proper morphism
of smooth quasi-projective varieties, and let $E \subset V$ be an
irreducible smooth divisor (proper over $S$) which meets
every fiber $V_s = \pi^{-1} (s)\quad (s \in S)$ transversally.
Let $E_s := V_s \cap E = \bigcup_{i=1}^mC_{s,i}$ be
the decomposition into irreducible components (ordered
arbitrarily for every $s \in S$).

\begin{enumerate}
\item[{\rm (a)}] Assume that
for every $s \in S$ and each $i=1, \ldots , m, \, 
C_{s,i} \simeq \pP^n \, (n\geq 1)$ with the conormal bundle
$J_{C_{s,i}}/J_{C_{s,i}}^2 \simeq \cO_{\pP^n}(1)$,
where $J_{C_{s,i}}$ is the ideal sheaf of 
$C_{s,i}$ in $V_s$.\footnote{In other words 
(by Kodaira's contraction theorem [Kod]) we
assume that each irreducible component $C_{s,i}$ can be
contracted in the fiber $V_s$ into a smooth point.}
Then there is a commutative diagram  

\begin{diagram}[w=7mm,h=8mm,midshaft]
V
&& \rTo^{\varphi} &&  W\\
&\rdTo<{\pi} &&
\ldTo>{\pi'}\\
&&  S \,
\end{diagram} 

\noindent where  
$W$ is a smooth
quasi-projective variety, $\pi '$ is a smooth morphism
and $\varphi$ is an $S$-contraction of the divisor
$E$ on $V$ onto a smooth subvariety
$A\subset W$ \'etale over $S$.

\item[{\rm (b)}] Let $E'$ be another smooth divisor on $V$ 
(proper over $S$)
which meets every fiber $V_s$ ($s \in S$) and the divisor $E$
transversally, and let $E_s':=E' \cap V_s= \bigcup_{j=1}^{m'}
C_{s,j}'$ be the decomposition into irreducible components. Assume
that for every $s \in S$ and for every pair $i,j$ such that
$C_{s,i} \cap C_{s,j}' \ne \emptyset$, this intersection
becomes a hyperplane under an isomorphism $C_{s,i} \simeq \pP^n$
as in (a). Then the restriction $\varphi |_{E'} :
E' \to W$ is an $S$-immersion onto a divisor (proper
over $S$) whose singularities are at worst transversal 
intersections (along $A$) of several smooth branches.
\end{enumerate}
\ethm

\no\proof (a) We may assume in the sequel that the base $S$ is connected.
Fix a very ample line bundle $H$ on $V$. For an arbitrary point $s_0\in S$,
letting $L=mH$ with $m$ sufficiently big, we may assume that
the line bundle $L|_{V_{s_0}} - K_{V_{s_0}}$ on the fiber 
$V_{s_0}$ over $s_0$ is ample, and so 
by the Kodaira vanishing theorem
\cite[Sect.1.2]{GriHar} we have 
$H^q(V_{s_0}, \cL_{s_0}) =0\quad\forall q\ge 1$
(here $\cL_{s}:={\cal O}_{V_s}(L)$). 
Since $\cL$ is a flat $\cO_S$-module, 
by the semi-continuity theorem
\cite[Thm. III.12.8]{Ha} it follows that 
\be\label{va} H^q(V_{s}, \cL_s) =0\quad\forall q\ge 1\ee
for every point $s$ in a 
neighborhood $S_0$ of the point $s_0$. 
Thus, for $m_0$ large enough and $L=L_0:=m_0 H$, (\ref{va}) holds 
for every point $s\in S$.

Since the divisor $E$ is irreducible, the monodromy of the
smooth family $\pi |_E : E \to S$ acts transitively on the
set of irreducible components $\{ C_{s,i} \}_{i=1}^m$ of
the fiber $E_s$ ($s \in S$). Hence all these components 
(regarded as cycles of $V$) are
algebraically (and then also numerically) equivalent.
Thus $k:= {\rm deg} (L_0|_{C_{s,i}})$ does not depend
on $s,i$. Consider the
(Cartier) divisors $L_j:=L_0+jE=m_0 H+jE$ on $V$ $(j \in \Z )$. 
Under our assumptions $\bigl(C_{s,i} \simeq \pP^n$ and $[C_{s,i}]|_{C_{s,i}}
\simeq \caO_{\pP^n} (-1)\bigr)$ for
every $s \in S$ and each $i=1, \dots , m$ we have
$L_j|_{C_{s,i}} \simeq \cO_{\pP^n}(l)$ with $l := k-j$, and so
\be\label{jva} H^q (E_{s}, (\cL_j)_s )=0 
\qquad \forall q\ge 1,\quad\forall j=0,\ldots,k\,.\ee 
Now the same argument as in the
proof of the Castelnuovo-Enriques-Kodaira theorem \cite[Appendix]{Kod}, 
\cite[p. 477]{GriHar} 
shows that: {\em

\smallskip\begin{enumerate}
\item[(i)] for every $j=0,\ldots,k$ 
and for every $s\in S$ the restriction map
$$ H^0 (V_{s}, (\cL_j)_s ) \to H^0(E_s, (\cL_j)_s) $$
is surjective,

\item[(ii)] the linear system $|L_k |_{V_s}|$
of divisors on $V_{s}$ is base point free, and

\item[(iii)] the associated morphism 
$\varphi_s : V_{s} \to \pP^{h-1}=\pP\bigl(H^0(V_s, (\cL_k)_s)^*\bigr)$ 
(with $h:=h^0(V_s,(\cL_k)_s)$) 
yields a contraction of the 
irreducible components $C_{s,i}$
$(i=1, \ldots , m)$ of the divisor 
$E_s\subset V_s$ into $m$ distinct smooth
points. 
\end{enumerate}}

\smallskip\no For the convenience of readers 
we sketch this argument.
For each $j=1,\ldots,k$ consider the short exact sequence of sheaves
$$0\to\cO_{V_s}(L_{j-1})\to  \cO_{V_s}(L_{j})
\to\cO_{E_s}(L_{j})\to 0$$
and the corresponding long exact cohomology sequence
\be\label{co} 0\to H^0(V_s,\,(\cL _{j-1})_s)\to  
H^0(V_s,\,(\cL _{j})_s)\to H^0(E_s,\,(\cL _{j})_s)\to 
H^1(V_s,\,(\cL _{j-1})_s)\to\ldots\,.\ee
It follows from (\ref{jva}) and (\ref{co}) 
that for every $q\ge 1$ the natural homomorphisms
$$H^q(V_s,\,(\cL _{0})_s)\to  H^q(V_s,\,(\cL _{1})_s)\to\ldots\to
H^q(V_s,\,(\cL _{k})_s)$$
are surjective, and so by (\ref{va}) 
all these groups vanish.
(In particular, 
\be\label{nva} H^q(V_s,
\,(\cL _{k-1})_s)=H^q(V_s,
\,(\cL _{k})_s)=0\qquad\forall 
q\ge 1,\quad\forall s\in S\,.{\rm )}\ee
Now (\ref{co}) implies (i). 

Since the line bundle $L_k|_{V_{s} \setminus E_s} 
\simeq L_0|_{V_{s} \setminus E_s}$
is very ample, the restriction 
$\varphi_s |_{V_{s} \setminus E_s}$ 
gives an embedding. Furthermore, 
the restriction $L_k|_{E_s}$  
is a trivial bundle, and so (ii) and (iii) easily follow. 

By Lemma \ref{l1}(c), (\ref{nva}) implies 
that the dimension $h:=h^0(V_s, (\cL_k)_s)$ 
is constant on $S$ and 
$\xi(L_k)=\bigl(\bigcup_{s\in S} H^0(V_{s}, (\cL_k)_{s})\to S\bigr)$
is a rank $h$ vector bundle on $S$. Thus, clearly, 
$$\varphi:V\to \pP(\xi(L_k)^*),\qquad \varphi\,|\,V_s=\varphi_s\,,$$
is a proper morphism onto a closed subvariety
$W:=\varphi (V)$ of (the total space of) the projective bundle
$\pP(\xi(L_k)^*)$; actually it 
consists of contracting 
the divisor $E\subset V$ onto a smooth subvariety
$A\subset W$ \'etale (and $m$-sheeted) over $S$ under the 
projection $\pi':={\rm pr}|_W$ where ${\rm pr}:\pP(\xi(L_k)^*)\to S$ 
is the standard projection. By \cite[Prop. II.7.10(b)]{Ha}
$\pP(\xi(L_k)^*)$ is a quasi-projective variety,
whence so is $W$. 
 
The morphism $\varphi: V\setminus E\to W\setminus 
A\hookrightarrow \pP^{h-1}$
is an embedding (indeed, so is the morphism given by the line bundle 
$L_k|_{V \setminus E} \simeq L_0|_{V \setminus E}$). 
Hence $ W\setminus A$
is a smooth variety. 
To show that the variety $W$ itself is smooth, we proceed locally
using local trivializations of the vector bundle $\xi(L_k)$. 
Fix a point 
$s \in S$ together with an affine neighborhood $S_0$ of $s$ in $S$,
and an index $i_0 \in
\{ 1, \ldots , m \}$. Since by (i) the restriction map
$H^0(V_s, (\cL_{k-1})_s) \to H^0 (E_s, (\cL_{k-1})_s)$ is
surjective, we can find $n+1$ sections
$\xi_{s,0} , \ldots , \xi_{s,n} \in H^0 (V_s, (\cL_{k-1})_s)$ 
which
are linearly independent when restricted to sections of 
$L_{k-1}|_{C_{s,i_0}} \simeq \cO_{\pr^n}(1)$. Fix also a section
$\eta_0$ of the line bundle $[E]$ over $V$ that
vanishes on $E$, another one $\eta_{s,1} \in H^0(V_s, (\cL_k)_s )$
which does not vanish on $C_{s,i_0}$, 
and a basis $\si_{s,j}\,\,\,(j=1,\ldots,h)$ 
of $H^0(V_s,(\cL_k)_s)$. Shrinking the neighborhood $S_0$ 
(if necessary) we may assume that $\si_{s,j}=\si_j|_{V_s}$ 
where $\si_j\in H^0(V_{S_0},\cL_k)\,\,\,(j=1,\ldots,h)$ 
(see Lemma \ref{l1}(b)) 
and for any point $s'\in S_0$
the restrictions 
$\si_j|_{V_{s'}}\in H^0(V_{s'},(\cL_k)_{s'})\,\,\,(j=1,\ldots,h)$ 
still form a basis.
Decomposing the sections 
$\eta_{s,1},\,\,(\eta_0|_{V_s})\cdot \xi_{s,l} \in H^0 
(V_s, (\cL_{k})_s)\,\,\,\,(l=0,\ldots,n)$ by
the basis 
$\si_{s,j}\,\,\,(j=1,\ldots,h)$
we may extend them to sections, say, 
$\eta_1,\,\xi_0,\ldots, \xi_n \in H^0(V_{S_0},\cL_k)$ decomposed 
by the system $\si_{j}\in H^0(V_{S_0},\cL_k)\,\,\,(j=1,\ldots,h)$ 
with the same coefficients. 
Then the ratios
\be\label{6} z_0':= {\frac {\xi_0}{\eta_1}}, \ldots ,
z_n':= {\frac {\xi_n}{\eta_1}} \, \ee 
can be pushed down to regular functions (say) $z_0,\ldots,z_n$
in a neighborhood of the point 
$c_{s,i_0} :=\varphi(C_{s,i_0})\in A$ in $W$ which give a 
local coordinate system on the fiber $W_s$
with center at the point 
$c_{s,i_0}$ (cf. \cite[Appendix]{Kod}, \cite[p. 477]{GriHar}). 
Clearly, they still give 
a local coordinate system on the fiber $W_{s'}$ 
around the point $c_{s',i_0} :=\varphi(C_{s',i_0})\in A$
close enough to $c_{s,i_0}$. Thus if $(x_1, \ldots , x_r)$
(with $r:= \dim_{\C} S$) is a local 
coordinate system at the point $s \in S$, then 
$(x_1, \ldots , x_r, \,z_0,\ldots ,z_n)$
define a local coordinate system on $W$
with center at the point 
$c_{s,i_0}$, and the projection $\pi'$ in these local coordinates is given as
$$(x_1, \ldots , x_r, \,z_0,\ldots ,z_n)\longmapsto (x_1, \ldots , x_r)\,.$$ 
This proves (a). 

\smallskip\no (b) Let $M$ be an algebraic vector bundle on an algebraic
variety $Z$, and let $L, L'$ be two transversal vector
subbundles of $M$, so that 
$M/(L \cap L') \simeq L/(L\cap L') \oplus L'/(L \cap L')$.
Then we have $M/L \simeq L'/(L \cap L')$. 
Letting $M:=TV|_{E \cap E'},  L:=TE|_{E \cap E'}$, and
$L' =TE'|_{E \cap E'}$, and using the above observation
we obtain an isomorphism of the normal bundles
\be\label{7} N_{(E\cap E')/E'} \simeq N_{E/V}.\ee
Since by our assumption, for every $s \in S$ and for every pair $(i,j)$
such that $C_{s,i} \cap C_{s,j}' \ne \emptyset$ we have
$(C_{s,i}, C_{s,i} \cap C_{s,j}') \simeq (\pP^n , \pP^{n-1})$ and
$N_{C_{s,i}/V_s} \simeq \cO_{\pP^n} (-1)$, it follows from (\ref{7})
that $N_{(C_{s,i} \cap C_{s,j}')/E_s'} \simeq \cO_{\pP^{n-1}} (-1)$.
This allows us to apply the argument as in (a) \footnote{which 
are also valid for $n=0$.}
replacing the pair $(V,E)$ by the pair
$(E', E \cap E')$. Hence the restriction $\varphi |_{E'}$ is an
$S_0$-contraction of the divisor $E \cap E'$ on $E'$. The 
image $\varphi ({E'})$ is a divisor in $W$ proper over $S_0$ and having
only smooth branches\footnote{In fact, for $n \geq 2$
the image itself is smooth since then
(assuming that $E'$ is smooth) for any $i \in \{ 1, \ldots, m \}$,
$C_{s,i} \cap C_{s,j}' \ne \emptyset$ for 
at most one value of $j$.}. 
Since for every $s \in S$ the fiber $E_s':= E' \cap V_s$ is a smooth 
divisor in $V_s$ and $\varphi_s=\varphi |_{V_s} : V_s \to W_s$ is the 
blowing up with the (finite) smooth center $A_s :=A \cap W_s$, it follows
that the intersection of local branches of the (reduced) divisor
$\varphi (E_s') \subset W_s$ at any point of $A_s$ is transversal.
Therefore, the branches of the divisor $\varphi (E') \subset W$
which contain the center $A$ meet also transversally each other
and every fiber $W_s \, (s \in S_0)$. Now the proof 
is completed.
\qed

\section{Combinatorial constructions}

We use below the following

\subsection{Terminology and notation}

\medskip Let $\pi:V\to S$ be a family of quasi-projective varieties.
{\sl Shrinking the base}
means passing to a new family $\pi|_{\pi^{-1}(U)} : \pi^{-1}(U)
\to U$ where $U$ is a Zariski open subset of $S$; 
usually we keep the same notation before and after shrinking the base.

By a {\em smooth family of quasi-projective varieties} we mean a smooth
surjective morphism $f : X \to S$ of smooth 
quasi-projective varieties; hereafter the base $S$ is supposed to be irreducible. 
Note that any quasi-projective family with a smooth total space
can be made smooth by shrinking the base.

We say that a family $\baf : V \to S$ is a {\em relative completion}
of $f:X\to S$ if $\baf$ is a proper morphism, $X \subset V$
is a Zariski open dense subset and $f= \baf |_X$. It is  {\em of simple
normal crossing} (or simply SNC) {\em type} if $D:= V \setminus X$
is a simple normal crossing divisor on $V$. 
If the family $\baf : V \to S$
is smooth
and each fiber $V_s:=\baf^{-1} (s), \,\,s \in S$, meets the divisor $D$
transversally along an SNC-divisor $D_s := D\cap V_s 
\subset V_s$, then we say that $(V,D)$ is a {\em relative SNC-completion}
of $X$. Clearly, any smooth relative completion with an SNC-divisor $D$
can be reduced to a relative SNC-completion by
shrinking the base.

Let $f: X \to S$  be a smooth family with all
fibers isomorphic to $\C^2$, and let $\baf : V \to S$ be its
relative SNC-completion. Then for every point $s \in S$ 
the `boundary divisor' $D_s$ is a {\sl rational
tree} (on the
smooth rational projective surface $V_s$). The latter means that each irreducible component 
$C_{s,i}$ of $D_s$ is a smooth rational 
curve, and the {\em weighted dual graph} (say) 
$\G_s$ of $D_s$ is
a tree (e. g., see \cite[\S2]{Zai}).

Let $v\in\G_s$ be an {\sl at most linear (-1)-vertex} 
of $\G_s$
(that is, the valence of $v$ is at most $2$ and the weight 
of $v$ is $-1$). 
The Castelnuovo contraction of the corresponding 
irreducible $(-1)$-component of $D_s$ leads again 
to an SNC-completion $(V_s',D_s')$ of $X_s\simeq\C^2$. 
The dual graph $\G_s'$ of $D_s'$ is obtained from 
$\G_s$ by the {\em blowing down} $v$.
The inverse operation on graphs is called a {\em blowing up}.
This blowing up (blowing down) is called {\sl inner}
if $v$ is a linear vertex of $\G_s$ and {\sl outer} if $v$ is terminal.
The graph $\G_s$ is called {\sl minimal} if no contraction is
possible; in this case it is linear \cite{Ram}. All minimal linear graphs 
corresponding to minimal SNC-completions of $\C^2$ 
are described in \cite{Ram}
and \cite{Mor}; we call them the
{\em Ramanujam-Morrow graphs}.

Assume that  
$\baf:(V,D)\to S$ as above is a proper and smooth 
SNC-family, and
fix a base point $s_0\in S$.  
There exists a smooth horizontal connection on $V$ 
which is tangent along the boundary SNC-divisor $D$
(indeed, it can be patched from local smooth connections 
tangent along $D$
using 
a smooth partition of unity on $V$). This provides us with 
a {\em geometric monodromy representation} 
$$\mu:\pi_1(S,\,s_0)\to {\rm Diff}\,(V_{s_0},\,D_{s_0})\,.$$
We denote by the same letter $\mu$ the induced 
{\em combinatorial monodromy representation} 
$ \pi_1(S,\,s_0)\to \Aut \G_{s_0}$. 

For a vertex $v\in \G_{s_0}$, let $O(v)$ be its $\mu$-orbit. 
Clearly, two vertices $v,\,v'\in\G_{s_0}$ belong to the same orbit 
if and only if,  for a certain irreducible component 
$E=E(v)$ of the boundary divisor $D\subset V$, the corresponding 
irreducible components $C(v)$ and $C(v')$
of the curve $D_{s_0}\subset V_{s_0}$ are contained 
in $E_{s_0}:=E\cap V_{s_0}$
(note that this fact is stable under shrinking the base). 
The next important lemma easily follows from Theorem \ref{pr1}.

\medskip\no\blem\label{bldown} In the notation as above, let $v$
be at most linear $(-1)$-vertex of the graph $\G_{s_0}$. 
Assume that the orbit $O (v) \subset \G_{s_0}$ of $v$ 
does not contain
a pair of neighbors in $\G_{s_0}$. 
Then (possibly, after  shrinking the base) 
there is a relative blowing down of the irreducible component 
$E(v)$ of the divisor $D$ which gives again 
a relative SNC-completion of the family $f:X\to S$.\elem

\subsection{Equivariant contractions}
From now on we consider a smooth SNC-completion 
of $\C^2$ by an SNC-divisor (say) $D_0$ with 
a weighted dual graph $\G$. 
Denote by $O(v)$ the orbit of a vertex 
$v$ of $\G$ under the action on $\G$ of 
the full automorphism group $\Aut \G$. 
If $v$ is at most linear $(-1)$-vertex such that its orbit 
$O(v)$ does not contain
a pair of neighbors in $\G$, then all the vertices in $O(v)$
can be simultaneously contracted; 
we call this an {\em equivariant contraction}
(or an {\em equivariant blowing down}). 
The main result of this subsection is the following proposition.

\medskip\no\bprop\label{pc} A graph $\G$ as above can be contracted to
a Ramanujam-Morrow graph by means of 
equivariant contractions.\eprop

The proof is done in Lemmas \ref{twolem} and \ref{pair} below.

\medskip\no\blem\label{twolem} Let $v_1$ and $v_2$ be  
at most linear $(-1)$-vertices of $\G$ which are neighbors
and belong to the same orbit (i.e., $v_2\in O(v_1)$). 
Draw $\G$ as follows:

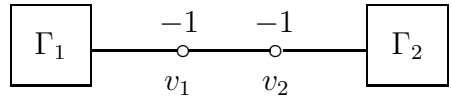
\begin{picture}(300,50)
\put(120,10){\framebox(30,30){$\G_1$}}
\put(151,23){\line(1,0){32}}
\put(182,20){$\circ$}
\put(178,8){$v_1$}
\put(176,30){$-1$}
\put(187,23){\line(1,0){31}}
\put(215,8){$v_2$}
\put(212,30){$-1$}
\put(217,20){$\circ$}
\put(223,23){\line(1,0){31}}
\put(255,10){\framebox(30,30){$\G_2$}}
\end{picture}

\noindent Then the following statements hold.

\smallskip \no (a) $\alpha(\G_i)=\G_j$ for any $\alpha\in\Aut \G$ with
$\alpha(v_i)=v_j,\,\,\,i,j\in\{1,2\}$. Moreover, $O(v_1)=\{v_1,v_2\}$,
and there is only one pair $\{v_1,v_2\}$ satisfying the assumptions of the lemma.

\smallskip \no (b) If $w\in\G_1$ is at most  linear $(-1)$-vertex
of $\G$ then the orbit $O(w)$ does not contain a pair of neighbors
in $\G$. \elem

\no\proof (a) Denote by ${\rm Br}_l(v_i)$ resp., ${\rm Br}_r(v_i)$
the left resp., the right-hand branch of $\G$ at $v_i$; thus
${\rm Br}_l(v_1)=\G_1$ and ${\rm Br}_r(v_2)=\G_2$.
Since ${\rm card}\,{\rm Br}_l(v_2)={\rm card}\,{\rm Br}_l(v_1)+1$,
for any $\alpha\in \Aut\G$ with $\alpha(v_1)=v_2$
we have $\alpha({\rm Br}_l(v_1))={\rm Br}_r(v_2)$, whence 
$\alpha(\G_1)=\G_2, \,\,\,\alpha(v_2)=v_1$ and 
$\alpha(\G_2)=\G_1$. 

In particular, ${\rm card}\,\G_1={\rm card}\,\G_2$. 
Therefore, if $\beta\in\Aut \G$ is such that $\beta(v_1)=v_1$ then 
$\beta(\G_1)=\G_1$, whence $\beta(v_2)=v_2$ and 
$\beta(\G_2)=\G_2$.
This proves the first statement of (a).

Let $(v_1',v_2')$ be another pair of at most 
linear $(-1)$-neighbors
of $\G$ which belong to the same orbit. As we have seen the both edges
$[v_1,v_2]$ and $[v_1',v_2']$ of $\G$ divide $\G$ into two parts of equal cardinality, which is only possible if $[v_1,v_2]=[v_1',v_2']$.
Thus the pair $(v_1,v_2)$ is unique. 
It follows that $O(v_1)=\{v_1,v_2\}$ which proves (a). 
It also follows that $O(w)\cap O(v_1)=\emptyset$ 
which proves (b).\qed
 
\no\blem\label{pair} Suppose that the graph $\G$ 
is not minimal.
Then $\G$ has an at most 
linear $(-1)$-vertex $w$ such that the orbit $O(w)$ 
does not contain a pair of neighbors.\elem

\no\proof Let $(v_1,v_2)$ be a pair of at most 
linear $(-1)$-neighbors
of $\G$ which belong to the same orbit. Clearly, $\G\neq \{v_1,v_2\}$, 
and so $\G$ contains a fragment

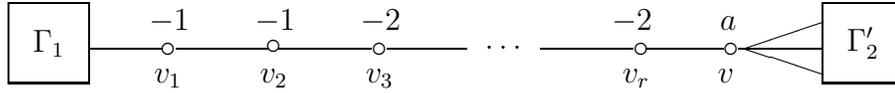
\begin{picture}(400,50)
\put(30,10){\framebox(30,30){$\G_1$}}
\put(60,23){\line(1,0){28}}
\put(87,20){$\circ$}
\put(83,30){$-1$}
\put(85,10){$v_1$}
\put(93,23){\line(1,0){35}}
\put(127,21){$\circ$}
\put(123,30){$-1$}
\put(125,10){$v_2$}  
\put(132,23){\line(1,0){35}}
\put(163,30){$-2$}
\put(167,20){$\circ$}
\put(165,10){$v_3$}
\put(173,23){\line(1,0){30}}
\put(210,23){$\ldots$}
\put(231,23){\line(1,0){35}}
\put(258,30){$-2$}
\put(266,20){$\circ$}
\put(262,10){$v_r$}
\put(271,23){\line(1,0){29}}
\put(299,30){$a$}
\put(300,20){$\circ$}
\put(298,10){$v$}
\put(306,23){\line(1,0){31}}
\put(308,23){\line(3,1){30}}
\put(308,23){\line(3,-1){30}}
\put(338,10){\framebox(30,30){$\G_2'$}}   
\end{picture}

\no where either $a\neq -2$ or $v$ is a branch vertex of $\G$.
Contracting the chain $(v_2,v_3,\ldots,v_r)$ 
we obtain the graph

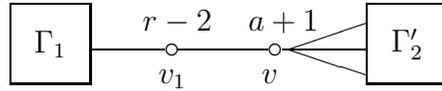
\begin{picture}(400,50)
\put(110,10){\framebox(30,30){$\G_1$}}
\put(140,23){\line(1,0){28}}
\put(168,20){$\circ$}
\put(160,30){$r-2$}
\put(166,10){$v_1$}  
\put(173,23){\line(1,0){34}}
\put(200,30){$a+1$}
\put(207,20){$\circ$}
\put(205,10){$v$}
\put(213,23){\line(1,0){30}}
\put(213,23){\line(1,0){31}}
\put(215,23){\line(3,1){30}}
\put(215,23){\line(3,-1){30}}
\put(245,10){\framebox(30,30){$\G_2'$}}   
\end{picture}

\no where the vertex $v$ cannot be contracted.
Assume that $\G$ does not contain at most linear
(-1)-vertices other than $v_1$ and $v_2$. 
Then after this contraction
the resulting graph is minimal
(that is, a Ramanujam-Morrow graph).
Show that this is impossible. Indeed, otherwise
$\G$ would be a linear graph admitting an automorphism 
$\alpha\in\Aut \G$ which interchanges $v_1$ and $v_2$ resp.,
$\G_1$ and $\G_2$. After the contraction as above, 
the image of $\G$ would
contain one of the following fragments
 
\begin{picture}(400,50)
\put(-10,20){$\cdots$}
\put(8,23){\line(1,0){7}}
\put(15,20){$\circ$}
\put(15,30){$a$}
\put(6,10){$\alpha(v)$}
\put(20,23){\line(1,0){7}}
\put(32,20){$\cdots$}
\put(50,23){\line(1,0){19}}
\put(69,20){$\circ$}
\put(60,30){$-2$}
\put(59,10){$\alpha(v_3)$}
\put(74,23){\line(1,0){21}}
\put(95,20){$\circ$}
\put(84,30){$r-2$}
\put(95,10){$v_1$}  
\put(100,23){\line(1,0){20}}
\put(116,30){$a+1$}
\put(120,20){$\circ$}
\put(120,10){$v$}
\put(125,23){\line(1,0){20}}
\put(152,20){$\cdots$} 
\put(172,20){$(r\ge 3),$}

\put(230,20){$\cdots$}
\put(250,23){\line(1,0){19}}
\put(269,20){$\circ$}
\put(268,30){$a$}
\put(260,10){$\alpha(v)$}
\put(274,23){\line(1,0){21}}
\put(295,20){$\circ$}
\put(295,30){$0$}
\put(294,10){$v_1$}  
\put(300,23){\line(1,0){20}}
\put(314,30){$a+1$}
\put(320,20){$\circ$}
\put(320,10){$v$}
\put(325,23){\line(1,0){20}}
\put(352,20){$\cdots$} 
\put(375,20){$(r=2).$}  
\end{picture}

\no  But a Ramanujam-Morrow graph can have only one positively weighted 
vertex, and a neighbor of this vertex has zero weight. Hence $a=-1$ 
in the left-hand fragment above, which contradicts the 
minimality assumption. 
Furthermore, the only fragments of a Ramanujam-Morrow graph of length
$3$ including a zero vertex in the middle are of the form

\begin{picture}(400,40)
\put(130,10){$\cdots$}
\put(150,13){\line(1,0){19}}
\put(169,10){$\circ$}
\put(169,20){$n$}
\put(174,13){\line(1,0){21}}
\put(195,10){$\circ$}
\put(195,20){$0$}  
\put(200,13){\line(1,0){20}}
\put(210,20){$-n-1$}
\put(220,10){$\circ$}
\put(225,13){\line(1,0){20}}
\put(252,10){$\cdots$}   
\end{picture}

\no where $n>0$ \cite[\S3.5]{Ram, Mor, FleZai}. 
Thus the both cases above are not possible.
In virtue of Lemma \ref{twolem} 
this concludes the proof. \qed

Proposition \ref{pc} and Lemma \ref{bldown} yield

\medskip\no\bcor\label{eq1} Let as above $f:X\to S$ 
be a smooth family of quasi-projective varieties with all fibers 
isomorphic to $\C^2$, and let $\baf:V\to S$ 
be its relative SNC-completion. Then the boundary 
divisor $D=V\setminus X$
(possibly, after shrinking the base) can be contracted 
providing a new  relative SNC-completion 
$\baf^{\rm min}:V^{\rm min}\to S$, where 
for each $s\in S$ the dual graph
$\G_s^{\rm min}$ of the boundary divisor 
$D_s^{\rm min}:=V_s^{\rm min}\setminus X_s$ is a 
Ramanujam-Morrow graph.\ecor

We need below the following lemma from \cite{FleZai}.

\medskip\no\blem \cite[L. 3.7]{FleZai}\label{re} 
Let $\G$ be a Ramanujam-Morrow graph. Then
$\G$ can be transformed, by a sequence of inner 
blowing ups and blowing downs, into one of the following
graphs:

\begin{picture}(400,40)
\put(50,10){$\circ$}
\put(50,20){$1$}

\put(129,10){$\circ$}
\put(130,20){$0$}
\put(134,13){\line(1,0){40}}
\put(173,20){$n$}
\put(173,10){$\circ$}   
\put(187,10){${\rm (} n\neq -1{\rm)}$}

\put(269,10){$\circ$}
\put(268,20){$0$}
\put(274,13){\line(1,0){21}}
\put(295,10){$\circ$}
\put(285,20){$k-1$}  
\put(300,13){\line(1,0){20}}
\put(319,20){$-1$}
\put(320,10){$\circ$}
\put(340,10){${\rm(}k\ge 1 {\rm)}\,.\qquad {\rm (*)}$}   

\end{picture}
 
\elem

\section{Proof of Theorem \ref{mt}} 
The next proposition is the key point 
in the proof of Theorem \ref{mt}.

\medskip\no\bprop\label{pr2} Let the assumptions 
of Corollary \ref{eq1} above be fulfilled. 
Then the family $f: X\to S$ 
(possibly, after shrinking the base)
admits a relative SNC-completion 
${\tilde f}:{\tilde V}\to S$ 
such that for every $s\in S$ 
the boundary divisor 
${\tilde D}_s:={\tilde V}_s\setminus X_s$ 
is irreducible
(and so isomorphic to $\pP^1$). \eprop 

\no\proof 
Let a relative SNC-completion
$(V^{\rm min},\,D^{\rm min})$ be as in Corollary \ref{eq1}.
Note that for any Ramanujam-Morrow graph $\G$
except the following one:

\begin{picture}(400,40)
\put(169,10){$\circ$}
\put(170,20){$0$}
\put(174,13){\line(1,0){40}}
\put(212,20){$0$}
\put(213,10){$\circ$}
\end{picture}

\no we have $\Aut \G=\{{\rm id}\}$. 
Let us deal with this exceptional case first. The
edge of this graph (invariant under  automorphisms) 
corresponds to a section (say) 
$\Sigma$ of $D^{\rm min}$ over $S$
($\Sigma$ is just the set of double points of the divisor $D^{\rm min}$).
We can blow up $V^{\rm min}$ along $\Sigma$ and then
(possibly, after shrinking the base)
blow down (according to Lemma \ref{bldown}) the proper transform(s)
of (the irreducible components of) $D$
to arrive at a new
relative SNC-completion with only irreducible
boundary divisors in fibers,
as required.

In the other cases the absence of nontrivial automorphisms implies
that 

\smallskip\no $\bullet$ each irreducible component $E_i$ of the boundary 
divisor $D^{\rm min}$ meets every fiber $V_s^{\rm min}$ along
an irreducible curve $C_{s,i}$, and 

\smallskip\no $\bullet$ the intersection of 
two such components 
$\Sigma_{ij}:=E_i\cap E_j \,\,\,(i\neq j)$
(if non-empty) is a (smooth) section. 

\smallskip\no These two properties are stable 
under blowing up with center at a section which is
the intersection of two components of the boundary divisor,
as well as under blowing down of a component of the boundary divisor
which corresponds to an at most linear (-1)-vertex
(it is defined correctly
over $S$ in virtue of  Lemma \ref{bldown}).

The last observation and Lemma \ref{re} imply that
a relative SNC-completion
$(V^{\rm min}$, $D^{\rm min})$ 
can be transformed
(possibly, after shrinking the base) into another
one with the dual graph
$\G_{s_0}$ as in $\rm (*)$.
If we finally arrive at a relative SNC-completion with 
the dual graph $\G=\G_{s_0}$ as in the third 
case of $\rm (*)$, then after blowing down 
the $(-1)$-curve in every fiber we obtain 
a relative SNC-completion with the dual graph
as in the second case of $\rm (*)$. 
In particular, we may assume that the singular 
locus $\Sigma$ of the boundary divisor 
$D$ is a section. 

If $n=0$ then we deal with the exceptional
case which is already settled.

If $n>0$ then we can proceed as before, 
performing first an inner 
relative blowing up with center at $\Sigma$ 
and then an outer relative blowing 
down. After a sequence of $n$ such `elementary 
transformations' we get
a relative SNC-completion with $n=0$, and so  we 
can finish the proof as above. 

Finally, consider the case where $n\le -2$. 
In this case we have $\Aut \G = \{{\rm id}\}$, and so the 
combinatorial monodromy of the family $\baf:D\to S$ is trivial.
Hence the divisor $D$ consists of two smooth 
irreducible components, say, $C_0$ and $C_1$ with $C^2_{s,0}=0$ 
and $C^2_{s,1}=n\le -2$ for every $s\in S$. 

Suppose that there exists a section $\Sigma'$ of 
$\baf\,|\,C_0:C_0\to S$
disjoint with $\Sigma:=C_0\cap C_1$. 
The kind of elementary transformations appropriate in our case 
is blowing up with center at $\Sigma'$ and then 
blowing down the proper transform of $C_0$ 
(by Theorem 
\ref{pr1}, this is possible after shrinking the base). 
Performing $n$ such elementary transformations 
(this needs at each step the existence of a section as above),
we arrive again at a relative SNC-completion 
of the second type with
$n=0$, and so we are done. 
Thus it remains to prove the following statement.

\medskip\no {\bf Claim.} {\em 
After shrinking the base appropriately
one can find a section $\Sigma'$ of 
$\baf\,|\,C_0:C_0\to S$
disjoint with $\Sigma:=C_0\cap C_1$.}

\medskip\no {\em Proof of the claim.} Letting in 
Lemma \ref{l1}
$E=C_0$ and $L=[C_1]|_E$ 
(so that $L|_{E_s}=L|_{C_{s,0}}\simeq\cO_{\pP^1}(1)$ 
for every $s\in S$)
and shrinking the base $S$ to make it affine, 
by Corollary \ref{cor1}(b)
we conclude that for every point $s\in S$ 
the restriction map
$$H^0(E, \cL) \to H^0(E_s, \cL_s )\simeq 
H^0(\pP^1,\cO_{\pP^1}(1))$$
is surjective. Thus for any points $s_0\in S$ and 
$z_0\in C_{s_0,0}\setminus C_{s_0,1}$ 
there exists a section 
$\sigma\in H^0(E, \cL)$ with 
$\si^*(0)\cdot C_{s_0,0}=z_0$.
The divisor $\Sigma':=\si^*(0)$ on $E=C_0$ 
(linearly equivalent to $\Sigma$)
passes through the point $z_0$ and meets 
every fiber $E_s=C_{s,0}\,\,\,(s\in S)$ 
transversally at one point. 
Clearly, $Z:=\baf(\Sigma\cap\Sigma')\not\ni s_0$ 
is a Zariski closed
proper subset of the base $S$. 
The restrictions of the sections 
$\Sigma$ and $\Sigma'$ onto the 
Zariski open subset $S_0:=S\setminus Z$
of $S$ are disjoint, as required. This proves the claim. \qed

Now the proof of Proposition \ref{pr2} is completed. 
\qed

In virtue of Corollary \ref{eq1} and Proposition \ref{pr2} 
the proof of Theorem \ref{mt} is reduced 
to the following simple lemma. It is well known that 
any smooth family with fibers isomorphic to a projective space
is locally trivial in the \'etale topology 
(and so it is a smooth Severi-Brauer variety)
\cite[Thm. I.8.2]{Gro}. It is locally trivial 
in the Zariski topology if and only 
if this family (or, equivalently, its dual) 
admits local sections; 
Lemma \ref{irr} below provides a proof along the lines 
in \cite[II, Sect. 0]{Gro}.

\medskip\no\blem\label{irr} Let $\baf:V\to S$ be 
a proper smooth family 
over a quasi-projective base with all fibers 
isomorphic to $\pP^n$, and let
$D$ be an irreducible smooth divisor on $V$ which 
meets every fiber $V_s\,\,\,(s\in S)$ transversally, 
with $D_s\simeq\pP^{n-1}$ and $N_{D_s/V_s}\simeq \cO_{\pP^n}(1)$.
Then the family $(V,D)$ is locally trivial in the Zariski topology.
\elem

\no\proof Fix an arbitrary point $s_0\in S$. 
Shrinking the base $S$ to an affine neighborhood 
of the point $s_0$, and letting in Lemma
\ref{l1} $E=V$ and $L=[D]$, by Corollary \ref{cor1} we will have that
the restriction map
$$H^0(E, \cL) \to H^0(E_{s_0}, \cL_{s_0})\simeq 
H^0(\pP^n,\cO_{\pP^n}(1))$$
is surjective. Fix sections $\si_0,,\ldots,\si_n\in H^0(E, \cL)$
such that their restrictions to the fiber 
$E_{s_0}$ are linearly independent and $\si_0^*(0)=D$. 
Shrinking the base further we may suppose that, 
for every fiber $E_{s}\,\,(s\in S)$, the restrictions 
$\si_0|_{E_s},\ldots,\si_n|_{E_s}$ are linearly independent as well. 
Then the morphism
$$\varphi:V\to S\times \pP^n,\qquad z\longmapsto 
\bigl(\baf(z),\,(\si_0(z):\ldots :\si_n(z))\bigr)$$
yields a desired trivialization over $S$.\qed

\end{document}